\documentclass[oneside]{amsart}

\usepackage[letterpaper,body={14.0cm,22.0cm}, mag=1000]{geometry}
\usepackage{amssymb}
\usepackage{amsthm}
\usepackage{amscd}
\usepackage{enumitem}
\usepackage{float}
\usepackage{placeins}
\usepackage{caption}
\numberwithin{equation}{section}
\theoremstyle{plain}

\newtheorem{thm}{Theorem}[section]
 \newtheorem{cor}[thm]{Corollary}
 \newtheorem{lemma}[thm]{Lemma}
\newtheorem{prop}[thm]{Proposition}

\theoremstyle{definition}

\newcommand{\dlabel}[1]{\ifmmode \text{\ttfamily \upshape [#1] } \else
{\ttfamily \upshape [#1] }\fi \label{#1}}

\newcommand{\Ho}{\operatorname{H} }
\newcommand{\M}{\operatorname{M} }

\newcommand{\Z}{\operatorname{Z} }

\newcommand{\Ker}{\operatorname{Ker} }

\newcommand{\im}{\operatorname{Im} }

\begin{document}

\title{Schur multipliers of special $p$-groups of rank $2$}

\author{Sumana Hatui}
\address{Department of Mathematics, Indian Institute of Science (IISc), Bangalore-560012, India}
\email{sumanahatui@iisc.ac.in, sumana.iitg@gmail.com}


\begin{abstract}
Let $G$ be a special $p$-group with center of order $p^2$. Berkovich and Janko asked to find the Schur multiplier of $G$ in \cite[Problem 2027]{B}. In this article we answer this question by explicitly computing the Schur multiplier of these groups. 
\end{abstract}

\maketitle 

\section{Introduction}
Let $G$ be a finite $p$-group and $d(G)$ denotes the cardinality of minimal generating set of $G$. The commutator subgroup and center of $G$ are denoted by $G'$ and $\Z(G)$ respectively. By $ES_{p^k}(p^n)$, we denote the extraspecial $p$-group of order $p^n$ of exponent $p^k$, $k=1,2$ and by $\mathbb{Z}_p^{(k)}$, we denote the elementary abelian $p$-group of rank $k$ for $ k \geq 1$.  For a prime $p$, $G^p$ denotes the group generated by the set $\{ g^p \mid g \in G \}$.
A finite $p$-group $G$ is called \emph{special $p$-group of rank $k$} if 
$G'=\Z(G)$ is an elementary abelian $p$-group of order $p^k$ and $G/G'$ is elementary abelian.
A group $G$ is called \emph{capable group} if there exists a group $H$ such that $G \cong H/\Z(H)$.
The \emph{epicenter} of a group $G$ is denoted by $Z^*(G)$, which is the smallest central subgroup $K$ of $G$ such that $G/K$ is capable.  

The Schur multiplier of a group $G$, denoted by $\M(G)$, is the second integral homology group $\Ho_2(G, \mathbb{Z})$ which was introduced by Schur in 1904 in his fundamental work on projective representation of groups.  There has been a great importance  in understanding the Schur multipliers of finite $p$-groups in recent past.  Here we are interested to   compute Schur multiplier of special $p$-groups of rank $2$. The special $p$-groups of minimum rank are the extraspecial $p$-groups and their Schur multiplier was studied in \cite{BE}. The Schur multiplier of special $p$-groups having maximum rank was studied in \cite{PK}. In this article we determine the Schur multiplier of special $p$-groups of rank $2$ and that answered the question which was asked by Berkovich and Janko in \cite[Problem 2027]{B}. 

Recall that, there are two extraspecial $p$-groups of order $8$ upto isomorphism, both are of exponent $4$, one is quaternian group $Q_8$ which has trivial Schur multiplier and another is dihedral group  $D_8$ which has Schur multiplier of order $2$ (see \cite[Theorem 3.3.6]{GK}).

We state our main results now. 
The following result describes the Schur multiplier  of $G$
when $G^p=G' \cong \mathbb Z_p \times \mathbb Z_p$.
\begin{thm}\label{S1}
Let $G$ be a special $p$-group of rank $2$ with $d=d(G)$ and $G^p=G'$. Then the following assertions hold:
\begin{enumerate}[label=(\alph*)]
\item Either $Z^*(G)=\Z(G)$ or $G$ is capable.

\item $\M(G)$ is elementary abelian of order $p^{\frac{1}{2}d(d-1)-2}$
if and only if $Z^*(G)=\Z(G)$. 

\item $\M(G)$ is of order $p^{\frac{1}{2}d(d-1)-1}$ of exponent at most $p^2$
if and only if $G$ is capable.
\item For every central subgroup $Z$ of order $p$, $G/Z$ is isomorphic to $$ES_{p^2}(p^{2m+1}) \times \mathbb{Z}_p^{(d-2m)}\;( m \geq 2), \; ES_{p^2}(p^3) \times \mathbb{Z}_p^{(d-2)}\; \text{or} \; Q_8 \times \mathbb{Z}_2^{(d-2)}.$$ 
\end{enumerate}
\end{thm}
By \cite[Proposition 3]{HH} and Theorem \ref{S1} we have the immediate corollary.
\begin{cor}
If $G$ is a special $p$-group of rank $2$ of order $p^n\; (n \geq 8)$ with $G^p=G'$, then
$\M(G)$ is elementary abelian of order $p^{\frac{1}{2}(n-2)(n-3)-2}$.
\end{cor}

The following result describes the Schur multiplier of $G$ when $G^p$ is cyclic  of prime order.
\begin{thm}\label{S2}
If $G$ is a special $p$-group of rank $2$ with $d=d(G)$ and $G^p \cong \mathbb{Z}_p$, then the following assertions hold:
\begin{enumerate}[label=(\alph*)]
\item $G$ is not capable and either $Z^*(G)=\Z(G)$ or $Z^*(G)=G^p$.


\item $\M(G)$ is elementary abelian.

\item The following are equivalent:
\begin{enumerate}[label=(\roman*)]
\item $\M(G)$ is of order $p^{\frac{1}{2}d(d-1)-2}$. 
\item $Z^*(G)=\Z(G)$.
\item $G/G^p \cong ES_p(p^{2m+1}) \times \mathbb{Z}_p^{(d-2m)}, m \geq 2$.
\end{enumerate}
\item The following are equivalent:
\begin{enumerate}[label=(\roman*)]
\item $\M(G)$ is of order $p^{\frac{1}{2}d(d-1)}$.
\item $Z^*(G)=G^p$.
\item $G/G^p \cong ES_p(p^{3}) \times \mathbb{Z}_p^{(d-2)}$.
\end{enumerate}
\end{enumerate}
\end{thm}
Now we are left with only one case when $G^p$ is trivial and in this case the  Schur multiplier of $G$  is studied in the following result.
\begin{thm}\label{S3}
If $G$ is a special $p$-group of rank $2$ with $d=d(G)$ and $G^p=1$, $p$ odd, then the following assertions hold:
\begin{enumerate}[label=(\alph*)]
\item $\M(G)$ is elementary abelian.

\item $p^{\frac{1}{2}d(d-1)-2} \leq |\M(G)| \leq p^{\frac{1}{2}d(d-1)+3}$.

\item $G$ is capable if and only if $G$ is isomorphic to one of the following groups: 

\begin{enumerate}[label=(\roman*)]
\item $\Phi_4(1^5)=\langle \alpha,\alpha_1,\alpha_2,\beta_1,\beta_2 \mid [\alpha_i,\alpha]=\beta_i, \alpha^p=\alpha_i^p=\beta_i^p=1 \; (i=1,2)\rangle$.

\item $\Phi_{12}(1^6)=ES_p(p^3) \times ES_p(p^3)$, 

\item $\Phi_{13}(1^6)=\langle \alpha_1,\alpha_2,\alpha_3,\alpha_4,\beta_1,\beta_2 \mid  [\alpha_i, \alpha_{i+1}]=\beta_i, [\alpha_2,\alpha_4]=\beta_2, \alpha_i^p=\alpha_3^p=\alpha_4^p=\beta_i^p=1 \;(i=1,2)\rangle$,

\item $\Phi_{15}(1^6)=\langle \alpha_1,\alpha_2,\alpha_3,\alpha_4,\beta_1,\beta_2 \mid  [\alpha_i, \alpha_{i+1}]=\beta_i, [\alpha_3,\alpha_4]=\beta_1,[\alpha_2,\alpha_4]=\beta_2^g, \alpha_i^p =\alpha_3^p = \alpha_4^p =\beta_i^p=1 \; (i=1,2)\rangle$,
where $g$ is non-quadratic residue modulo $p$.

\item $T =\langle x_1,\cdots , x_5, c_1, c_2 \mid [x_2, x_1] =[x_5, x_3]= c_1, [x_3, x_1] = [x_5, x_4] = c_2, x_i^p=c_j^p=1, \; 1 \leq i \leq 5, 1 \leq j \leq 2 \rangle$.
\end{enumerate}

\item $|\M(G)| = p^{\frac{1}{2}d(d-1)+3}$ if and only if $G$ is isomorphic to $\Phi_4(1^5)$.

\item $|\M(G)| = p^{\frac{1}{2}d(d-1)+2}$ if and only if $G$ is isomorphic to $\Phi_{12}(1^6),  \Phi_{13}(1^6), \; \text{or}\; \Phi_{15}(1^6)$.

\item  $|\M(G)|=p^{\frac{1}{2}d(d-1)-1}$ if and only if $G$ is isomorphic to $T$.


\item $\M(G)$ is of order $p^{\frac{1}{2}d(d-1)-2}$ if and only if $Z^*(G)=\Z(G)$. In this case, $G/Z \cong ES_p(p^{2m+1}) \times \mathbb{Z}_p^{(d-2m)}, m \geq 2$, for every central subgroup $Z$ of order $p$.


\item $\M(G)$ is of order $p^{\frac{1}{2}d(d-1)}$ if and only if  $Z^*(G)\cong\mathbb{Z}_p$. In this case, $G/Z^*(G) \cong ES_p(p^{3}) \times \mathbb{Z}_p^{(d-2)}$.
\end{enumerate}
\end{thm}
The following result is for $p=2$.
\begin{thm}\label{p2}
Let $G$ be a special $2$-group of rank $2$. Then $G^2=G'$  holds.
\end{thm}
\section{Preliminaries}\label{Ra}
For a finite group $G$ of class $2$ with  $G/G'$ elementary abelian, the following construction is given in \cite{BE}.  We consider $G/G'$ and $G'$ as vector spaces over $\mathbb{F}_p$, which we denote by $V, W$ respectively. The bilinear map $(-,-) : V \times V \to W$ is defined by 
$$(v_1,v_2)=[g_1,g_2]$$
for  $v_1,v_2 \in V$ such that $v_i=g_iG', i \in \{1,2\}$ for some $g_1,g_2 \in G$.  
Let $X_1$ be the subspace of $V \otimes W$ spanned by all 
$$v_1 \otimes (v_2,v_3) + v_2 \otimes (v_3,v_1) + v_3 \otimes (v_1,v_2)$$
for $v_1, v_2, v_3 \in V$.
Consider a map $f : V \rightarrow W$ given by 
$$f(gG')=g^p$$
 for $g \in G$.
Let $X_2$ be the subspace spanned by all $v \otimes f(v)$, $v \in V$,  and take $$X:=X_1+X_2.$$
Now consider a homomorpism $\sigma: V\wedge V \rightarrow (V \otimes W)/X$ given by
$$\sigma(v_1 \wedge v_2)= \big(v_1 \otimes f(v_2)+ (_2^p)  v_2 \otimes (v_1,v_2)\big)+X.$$
Then there exists an abelian group $M^*$  with a subgroup $N$ isomorphic to  $(V \otimes W)/X$, such that 
$$1 \rightarrow N \rightarrow M^* \xrightarrow{\xi} V \wedge V \rightarrow 1$$
is exact and $$\sigma\xi(\alpha)=\alpha^p \;\text{for}\; \alpha \in M^*.$$
Now we consider a homomorphism $\rho:V \wedge V \rightarrow W$ given by
$$\rho(v_1 \wedge v_2)=(v_1,v_2)$$ 
for all $v_1,v_2 \in V.$  Notice that $\rho$ is an epimorphism.
We let $M$ be the subgroup of $M^*$ containing $N$ such that $M/N \cong \Ker \rho$. We use this notation throughout the paper without further reference.

With the above setting, we have 
\begin{thm}$($\cite[Theorem 3.1]{BE}$)$ \label{thmBE}
 $\M(G) \cong M$.
\end{thm}
\textbf{Note:}
It is easy to observe that $X_1$ is generated by the set $$\{\bar{x}_1 \otimes [x_2,x_3]+ \bar{x}_2 \otimes [x_3,x_1]+ \bar{x}_3 \otimes [x_1.x_2] \mid \;x_1,x_2,x_3 \in S\},$$
where $S$ is a set of generators of $G$ and $\bar{x}$ is the image of $x$ in $G/G'$.

Suppose $G$ has a free presentation $F/R$. Let $Z=S/R$ be a central subgroup of $G$. 
Then the map from $(F/F'R) \times (S/R)$ to $(F'\cap R)/[F,R]$ defined by
 $(xF'R,sR) \mapsto [x,s][F,R]$ is a well-defined bilinear map and induces a homomorphism $\lambda_Z: (G/G') \otimes Z \to \M(G)$, called the Ganea map.
\begin{thm}[\cite{GT}]\label{G}
Let $Z$ be a central subgroup of a finite group $G$. Then the following sequence is exact 
$$ (G/G') \otimes Z \xrightarrow{\lambda_Z} \M(G) \xrightarrow{\mu} \M(G/Z) \rightarrow G'\cap Z \rightarrow 1.$$
\end{thm}
\begin{thm}[\cite{F}]\label{G1}
Let $Z$ be a central subgroup of a finite group $G$. Then $Z \subseteq Z^*(G)$ if and only if $(G/G') \otimes Z=\Ker \lambda_Z$.
\end{thm}
By \cite[Corollary 3.2.4]{GK}, we have $X = \Ker \lambda_{\Z(G)}$. Hence by Theorem \ref{G1}, we have the following result:
\begin{lemma}\label{Ka}
Let $Z$ be a central subgroup of a group $G$ of nilpotency class $2$. Then $Z \subseteq Z^*(G)$ if and only if $(G/G') \otimes Z$ is contained in $X$. 
\end{lemma}
Let $v_1, v_2, \ldots, v_d$ be the generators of $V$ such that $\{f(v_1),f(v_2), \ldots, f(v_r)\}$ is a basis of $G^p$. Then the set 
 $$\{v_i \otimes f(v_i), v_l \otimes f(v_i), (v_i \otimes f(v_j)+v_j \otimes f(v_i)) \mid 1 \leq i < j \leq r, (r+1) \leq l \leq d\}$$ forms a basis of $X_2$, from which the following result follows.
\begin{prop}[Proposition 3.3 of \cite{PK}]\label{PK}
Let $G$ be a special $p$-group with $d=d(G)$ and $G^p$ of order $r$. Then $|X_2|=p^{rd-\frac{1}{2}r(r-1)}$.
\end{prop}
By Theorem \ref{G}, we have 
\begin{eqnarray}\label{3}
\frac{|\M(G)|}{|\im \lambda_Z|}=\frac{|\M(G/Z)|}{|G' \cap Z|}.
\end{eqnarray}
As $X =\Ker \lambda_{\Z(G)}$, so $|\im \lambda_{\Z(G)}|=\frac{p^{2d}}{|X|}$.
Hence, by (\ref{3}), taking $Z=\Z(G)$, we have
\begin{eqnarray}\label{1}
|\M(G)|=p^{\frac{1}{2}d(d-1)-2}. \frac{p^{2d}}{|X|}
\end{eqnarray}
Now we recall the following results which will be used in the proof of the main results.
\begin{thm}[Main Theorem of \cite{PN}]\label{N3}
Let $G$ be a $p$-group of order $p^n$. Then $|\M(G)|=p^{\frac{1}{2}(n-1)(n-2)+1}$ if and only if $G \cong ES_p(p^3) \times \mathbb{Z}_p^{(n-3)}$.
\end{thm}
\begin{thm}[Theorem 21 of \cite{PN2}]\label{N1}
Let $G$ be a $p$-group of order $p^n$. Then $|\M(G)|=p^{\frac{1}{2}(n-1)(n-2)}$ if and only if $G$ is isomorphic to one of the following groups.
\begin{enumerate}[label=(\roman*)]
\item $G\cong \mathbb{Z}_{p^2} \times \mathbb{Z}_p^{(n-2)}$,

\item $G\cong  D_8 \times \mathbb{Z}_2^{(n-3)}$,

\item $G\cong  \mathbb{Z}_p^{(4)} \rtimes \mathbb{Z}_p\; (p \neq 2)$.
\end{enumerate}
\end{thm}
\begin{thm}[Theorem 11 of \cite{PN4}]\label{N2}
Let $G$ be a group of order $p^n$ such that $G/G'$ is elementary abelian of order $p^{n-1}$. Then $|\M(G)|=p^{\frac{1}{2}(n-1)(n-2)-1}$ if and only if $G$ is isomorphic to one of the following groups.
\begin{enumerate}[label=(\roman*)]
\item $ES_{p^2}(p^3) \times \mathbb{Z}_p^{(n-3)}$,

\item $Q_8 \times \mathbb{Z}_2^{(n-3)}$,

\item $ES(p^{2m+1}) \times \mathbb{Z}_p^{(n-2m-1)}\; (m \geq 2)$.
\end{enumerate}
\end{thm}
\section{Proofs}
In this section we prove our main results.

\vspace{.2in}

\noindent {\bf\emph{Proof of Theorem \ref{S1}:}}
Consider $G^p=G' \cong \mathbb{Z}_p \times \mathbb{Z}_p$. 
Now by Proposition \ref{PK}, $|X_2|=p^{2d-1}$. 
Let $v_1,v_2 \in V$ such that $\{f(v_1),f(v_2)\}$ is a basis of $W$.
Observe that $X_2$ is generated by the set
$$\{v_i \otimes f(v_i), v_j \otimes f(v_i),\;  \big(v_1 \otimes f(v_2)+v_2 \otimes f(v_1)\big) \mid  i=1,2 \; \text{and} \; 3 \leq j \leq d\}.$$ Hence $ p^{2d}\geq |X| \geq p^{2d-1}$ and by (\ref{1}),
$$p^{\frac{1}{2}d(d-1)-2} \leq |\M(G)| \leq p^{\frac{1}{2}d(d-1)-1}.$$

$(a)$ Observe that $|X|=p^{2d}$, i.e., the set $\{v_1 \otimes f(v_2), v_2 \otimes f(v_1)\}$ is contained in $X$ if and only if $Z^*(G)=\Z(G)$,  follows by Lemma \ref{Ka}. 
Another possibility  is $|X|=p^{2d-1}$, i.e., $X_1 \subseteq X_2$ if and only if $v_1 \otimes f(v_2), v_2 \otimes f(v_1)$ are not in $X$. Hence, by Lemma \ref{Ka}, $Z^*(G)=1$. So $G$ is capable.

$(b)$ By (\ref{1}), $|X|=p^{2d}$ if and only if
$$|\M(G)|=p^{\frac{1}{2}d(d-1)-2}.$$ By Theorem \ref{G1} and Theorem \ref{G}, taking $Z=\Z(G)$, we see that $\M(G)$ embeds in $\M(G/\Z(G))$ which is elementary abelian. So in this case $\M(G)$ is elementary abelian. 

$(c)$ For $x \in \M(G), \; x^p \in (V \otimes W)/X$ and $V \otimes W/X$ is elementary abelian, so $x^{p^2}=1$.

By (\ref{1}),  $X_1 \subseteq X_2$ i.e., $|X|=p^{2d-1}$ if and only if, by (\ref{1}), 
$$|\M(G)|=p^{\frac{1}{2}d(d-1)-1}.$$ Hence, by Lemma \ref{Ka}, $Z^*(G)=1$. So $G$ is capable. Converse follows from $(b)$.




$(d)$ Observe that, in both the cases $|\M(G/Z)|=p^{\frac{1}{2}d(d-1)-1}$, follows from (\ref{3}), taking $Z$ a central subgroup of order $p$. Therefore by Theorem \ref{N2},
\begin{eqnarray*}
G/Z &\cong& ES_{p^2}(p^{2m+1}) \times \mathbb{Z}_p^{(d-2m)}, \;(m \geq 2), ES_{p^2}(p^{3}) \times \mathbb{Z}_p^{(d-2)} \; \text{or}\; Q_8 \times \mathbb{Z}_2^{(d-2)}.
\end{eqnarray*} 
\hfill $\Box$

\vspace{.2 in}

\noindent {\bf\emph{Proof of Theorem \ref{S2}:}}
$(a)$ Assume $G^p$ is cyclic of order $p$. 
Now by Proposition \ref{PK}, $|X_2|=p^d$. 
Let $v_1\in V$ such that $G^p=\langle f(v_1) \rangle$.
Observe that $X_2$ is generated by the set
$$M':=\{v_i \otimes f(v_1), \; 1 \le i\leq d\}.$$ Hence by Lemma \ref{Ka}, $G^p\subseteq Z^*(G)$. 
Hence $G$ is not capable with $Z^*(G)=G^p$ or $Z^*(G)=\Z(G)$. 

Using Theorem \ref{G1} and taking $Z=G^p$ in (\ref{3}), we have 
\begin{eqnarray}\label{2}
|\M(G)|=\frac{|\M(G/G^p)|}{p}
\end{eqnarray}
Now let $v_1,v_2 \in V$ such that $(v_1,v_2) \in G'\setminus G^p$. Then The set 
$$N':=\{v_i\otimes(v_1,v_2)+v_1\otimes(v_2,v_i)+v_2\otimes(v_i,v_1)\mid 3\leq i \leq d\}$$ 
is linearly independent in $X_1$, so $|X_1| \geq p^{d-2}$. 
The set $M' \cup N'$ is linearly independent in $X$ and $M' \cap N'=\emptyset$.
Thus $$p^{2d} \geq |X| \geq p^{2d-2}.$$ 
Hence by (\ref{1}), 
$$p^{\frac{1}{2}d(d-1)-2} \leq |\M(G)|\leq p^{\frac{1}{2}d(d-1)}.$$

$(c)$ Now similarly, as described in the proof of Theorem \ref{S1}, we have $Z^*(G)=\Z(G)$ if and only if  $|\M(G)|=p^{\frac{1}{2}d(d-1)-2}$ i.e., by (\ref{2}), 
$$|\M(G/G^p)|=p^{\frac{1}{2}d(d-1)-1},$$ which happens 
if and only if $$G/G^p \cong ES_p(p^{2m+1}) \times \mathbb{Z}_p^{(d-2m)}, m \geq 2,$$
follows from Theorem \ref{N2}.

$(d)$ By Theorem \ref{N1}, it follows that there is no $G/G^p$ such that 
$$|\M(G/G^p)|=p^{\frac{1}{2}d(d-1)}.$$ Thus by (\ref{2}),
$|\M(G)|$ cannot be of order $p^{\frac{1}{2}d(d-1)-1}$.
Hence $|\M(G)|=p^{\frac{1}{2}d(d-1)}$ if and only if $Z^*(G)=G^p$.  By (\ref{2}), $$|\M(G/G^p)|=p^{\frac{1}{2}d(d-1)+1},$$
which happens if and only if
$$G/G^p \cong ES_p(p^{3}) \times \mathbb{Z}_p^{(d-2)},$$
follows by Theorem \ref{N3}.

$(b)$  
By $(c)$ and $(d)$ it follows that $p$ must be odd.
The group $G/G^p$ is of exponent $p$ and $p$ odd, so 
the homomorphism $\sigma$, described in Section \ref{Ra}, is trivial map and therefore, $\sigma \xi(x)=x^p=1$ for $x \in \M(G/G^p)$. Thus $\M(G/G^p)$ is elementary abelian. Since $G^p \subseteq Z^*(G)$, by Theorem \ref{G1} and Theorem \ref{G}, $\M(G)$ embeds in $\M(G/G^p)$. Therefore, $\M(G)$ is also elementary abelian.
The proof is complete now.
\hfill $\Box$

\vspace{.2 in}

\noindent {\bf\emph{Proof of Theorem \ref{S3}:}}
$(a)$ Since $p$ is odd and $G^p=1$, the homomorphism $\sigma$ described in Section \ref{Ra}, is the trivial map and therefore $\sigma \xi(x)=x^p=1$.
Thus, $\M(G)$ is elementary abelian.

Let $z,z'$ be the generators of $\Z(G)$ and $x_1,x_2, \ldots, x_d$ be the generators of $G$ such that $[x_1,x_2] \in \langle z \rangle$ is non-trivial.
Then the set
$$A:=\{x_i \otimes [x_1,x_2]+x_1 \otimes[x_2,x_i]+x_2 \otimes [x_i,x_1]\mid 3 \le i \le d\}$$
consists of $d-2$ linearly independent elements of $X_1$.

Now if for some $x_k  \in \{3,4, \ldots, d\}$, $[x_1,x_k] \in \langle z' \rangle$ is non-trivial, then the set
$$B :=\{x_i \otimes [x_1,x_k]+x_1 \otimes[x_k,x_i]+x_k \otimes [x_i,x_1]\mid 3 \le i \le d, i \neq k\}$$
consists of $d-3$ linearly independent elements of $X_1$. Thus, $A \cup B$ consists of $(2d-5)$ linearly independent elements of $X_1$. 

If $1 \neq [x_2,x_k]\in \langle z' \rangle$, for some $k \in \{3,4, \ldots, d\}$, then a similar conclusion holds. 
Suppose then that $[x_1,x_k]$ and $[x_2,x_k]$ are all trivial or in $\langle z \rangle$ for all $k  \in \{3,4, \ldots, d\}$. Say, $[x_3,x_4] \in \langle z' \rangle$.
In this case
$$B_1:=\{x_i \otimes [x_3,x_4]+x_3 \otimes[x_4,x_i]+x_4 \otimes [x_i,x_3]\mid 1 \le i \le d, \; i \neq 3,4 \}$$
consists of $d-2$ independent elements of $X_1$. Thus, $A \cup B_1$ consists of $2d-4$ linearly independent elements of $X_1$.

Hence in both cases, $2d \ge |X|=|X_1| \geq p^{2d-5}$ holds.
By (\ref{1}), it follows that,
$$p^{\frac{1}{2}d(d-1)-2} \leq |\M(G)| \leq p^{\frac{1}{2}d(d-1)+3},$$ which proves $(b)$.

Now if $G$ is not capable, then by (\ref{3}) and \cite[Main Theorem]{PN} we have  $$p^{\frac{1}{2}d(d-1)-2} \leq |\M(G)| = \frac{|\M(G/Z)|}{p} \leq p^{\frac{1}{2}d(d-1)}\; \text{for}\; Z \subseteq Z^*(G),\; Z \cong \mathbb{Z}_p.$$ 
By Theorem \ref{N1}, there is no $G$ and central subgroup $Z$ such that  $|\M(G/Z)|=p^{\frac{1}{2}d(d-1)}$.
Hence 
\begin{eqnarray}\label{4}
|\M(G)|=p^{\frac{1}{2}d(d-1)-2} \; \text{or} \; p^{\frac{1}{2}d(d-1)}, \;\text{if $G$ is not capable.}  
\end{eqnarray}

Assume then that $G$ is capable. By \cite[Proposition 3]{HH}, $p^5 \leq |G| \leq p^7,$ with $p$ odd. 
If $|G|=p^5$, then looking through the list of groups given in \cite{RJ}, it follows that $G \cong \Phi_4(1^5)$. Since $|X|=p$, it follows by (\ref{1}) that $|\M(G)|=p^6=p^{\frac{1}{2}d(d-1)+3}$.

If $|G|=p^6$, then looking through the list of groups given in \cite{RJ}, it follows that $G \cong \Phi_{12}(1^6), \; \Phi_{13}(1^6)$ or $\Phi_{15}(1^6)$. Since $|X|=p^4$, it follows by (\ref{1}) that $|\M(G)|=p^8=p^{\frac{1}{2}d(d-1)+2}$.

Now consider groups of order $p^7$ of exponent $p$. By \cite{HKM},  it follows that there is only one capable group 
\begin{eqnarray*}
G =\langle x_1,\ldots , x_5, c_1, c_2 \mid [x_2, x_1] =[x_5, x_3]= c_1, [x_3, x_1] = [x_5, x_4] = c_2, x_i^p=c_j^p=1, \\
1 \leq i \leq 5, 1 \leq j \leq 2 \rangle
\end{eqnarray*}
 up to isomorphism.
By (\ref{1}), $|\M(G)|=p^9=p^{\frac{1}{2}d(d-1)-1}$.
 
Now $(c), (d), (e), (f)$ follow by (\ref{4}).

$(g)$ Now by (\ref{3}) we have that $Z^*(G) \cong \Z(G)$ if and only if $|\M(G)|=p^{\frac{1}{2}d(d-1)-2}$,  taking $Z=\Z(G)$.
Hence by (\ref{3}) we have that for every central subgroup $Z$ of order $p$, $|\M(G/Z)|=p^{\frac{1}{2}d(d-1)-1}$, thus  $$G/Z \cong ES_p(p^{2m+1}) \times \mathbb{Z}_p^{(d-2m)},\; m \geq 2,$$ follows from Theorem \ref{N2}.

$(h)$ Suppose $Z^*(G) \cong \mathbb{Z}_p$. 
By (\ref{4}) and $(g)$, it follows that $|\M(G)|=p^{\frac{1}{2}d(d-1)}$.

Conversely suppose $|\M(G)|=p^{\frac{1}{2}d(d-1)}$.
From the previous cases it follows that $G$ is not capable since $Z^*(G) \cong \mathbb{Z}_p$.

In this case, we have $$|\M(G)| = \frac{|\M(G/Z^*(G))|}{p}.$$
By Theorem \ref{N3}, it follows that $|\M(G/Z^*(G))|=p^{\frac{1}{2}d(d-1)+1}$ if and only if 
$$G/Z^*(G) \cong ES_p(p^{3}) \times \mathbb{Z}_p^{(d-2)}.$$
The proof is complete now.
\hfill $\Box$

\vspace{.2in}
\noindent {\bf\emph{Proof of Theorem \ref{p2}:}}
By Theorem \ref{S2}, it follows that there is no special $2$-group of rank $2$ with $G^2 \cong \mathbb{Z}_2$.

Assume that $p=2$. As in the proof of Theorem \ref{S3}, we conclude that
$$p^{\frac{1}{2}d(d-1)-2} \leq |\M(G)| \leq p^{\frac{1}{2}d(d-1)+3}.$$
Let $Z$  be a central subgroup of order $2$.
If $G$ is not capable then by (\ref{3}),  
$$2^{\frac{1}{2}d(d-1)-2} \leq |\M(G)| = \frac{|\M(G/Z)|}{2} \leq 2^{\frac{1}{2}d(d-1)}\; \text{for}\; Z \subseteq Z^*(G).$$ 
By Theorem \ref{N3}, \ref{N1} and  \ref{N2}, there is no group $G$ and central subgroup $Z$ such that  $G/Z$ of exponent $2$ and $2^{\frac{1}{2}d(d-1)-1} \leq |\M(G/Z)| \leq 2^{\frac{1}{2}d(d-1)+1}$.
Hence $G$ must be capable and $|\M(G/Z)| \leq 2^{\frac{1}{2}d(d-1)-2}$.

Suppose $|\M(G)|=2^{\frac{1}{2}d(d-1)-2+k}, \; 0 \leq k \leq 5$. Then by (\ref{1}), $|\Ker \lambda_{\Z(G)}|=|X|=2^{2d-k}$. Hence 
$1 \leq |\im \lambda_Z|=2^m \leq 2^k$. By (\ref{3}), it follows that 
$$2^{\frac{1}{2}d(d-1)-2+k-m}=\frac{|\M(G)|}{|\im \lambda_Z|} = \frac{|\M(G/Z)|}{2} \leq 2^{\frac{1}{2}d(d-1)-3}.$$ Hence $k-m \leq -1$, which is not possible. So there is no special $2$-group of rank $2$ with $G^2=1$. 
The proof is complete now.

\hfill$\Box$

\textbf{Acknowledgement:} The research of the author is partly supported by Infosys grant.

\end{document}